\documentclass[11pt]{article}
\usepackage{graphicx}

\newcounter{Prop}[section]
\newcounter{Define}[section]
\renewcommand{\theProp}{\thesection.\arabic{Prop}.}
\renewcommand{\theDefine}{\thesection.\arabic{Define}.}

\newenvironment{Theorem}[1]
{\refstepcounter{Prop} \\[5pt] {\bf Theorem \theProp\ (#1)}\ \ }
{\\[5pt]}
\newenvironment{Proposition}
{\refstepcounter{Prop} \\[5pt] {\bf Proposition \theProp}\ \ }
{\\[5pt]}
\newenvironment{Corroraly}
{\refstepcounter{Prop} \\[5pt] {\bf Corroraly \theProp}\ \ }
{\\[5pt]}
\newenvironment{Definition}[1]
{\refstepcounter{Define} \\[5pt] {\bf Definition \theDefine\ (#1)}\ \ }
{\\[5pt]}

\title{A sufficient condition for pseudointegrable systems with weak mixing property}
\author{Minoru Ogawa\thanks{electric address: mogawa@r.phys.nagoya-u.ac.jp} \\
Department of Physics, Nagoya University, \\ Nagoya 464-8602 JAPAN}

\date{}

\begin{document}
\maketitle

\begin{abstract}
We present a sufficient condition that a pseudointegrable system has weak mixing property.
The result is derived from Veech's weak mixing theorem for interval exhange
[Veech, W.A. Amer.J.Math. {\bf 106}, 1331 (1984)].
We also present an example whose weak mixing property can be proved by the result.
\end{abstract}

\section{Introduction}
Billiards have attracted much attention in the recent decades as simple dynamical systems
in classical and quantum mechanics. Depending on shapes of boundary, statistical property
of billiards varies from integrable to strongly chaotic. \par
Plane polygonal billiards, or simply, polygonal billiards, is a typical class that is
neither chaotic nor integrable except to a few integrable cases (regular triangular,
($\pi /6,\ \pi /3,\ \pi /2$)-triangular, ($\pi /4,\ \pi /4,\ \pi /2$)-triangular and
rectangular tables only make their billiards integrable).
Dynamics of billiards in a typical polygon is conjectured to be ergodic on the three
dimensional energy surfaces. While dynamics of billiards in rational polygons are
restricted to two dimensional invariant surfaces as same as integrable systems, but genuses
of the invariant surfaces are larger than 1 (except to integrable cases).
Therefore billiards in rational polygons are characterized as pseudointegrable \cite{RB}.
It is rigorously proven that the dynamics on these surfaces are ergodic and not mixing
\cite{G}. \par
Every pseudointegrable system has continuous component of spectrum and may (or may not)
have {\it non-trivial} discrete one because of its non-integrable and non-chaotic properties.
Here a non-trivial discrete component of spectrum means a discrete one
which does not correspond to almost everywhere constant functions.
It turns out that a dynamics without non-trivial discrete component of spectrum has weak
mixing property. Here weak mixing property for a finite ergodic measure $\mu $ is defined as
follows: If the condition
\[
\lim_{T\to\infty }\frac{1}{T}\int_0^T \left| (U^tf,g)-(f,1)(1,g)\right| \mathrm{d}t=0
\]
is held for any $\mu $-square integrable functions $f,g$,
then it is said that the system has (measure theoretical) weak mixing property.
Here $(f,g)$ denotes $\int f\cdot g^\ast\mathrm{d}\mu $,
which represents the inner product on the Hilbert space composed of $\mu$-square integrable
functions.
A weak mixing dynamics without stronger mixing has singular continuous spectrum. \par
Weak mixing is believed to be generic for billiards in rational polygons,
or more generally, for geodesic flows on translation surfaces which is explained in next
section, but this was only shown for {\em horizontal-vertical polygons} or
{\em horizontal-vertical translation surfaces}. (See the Appendix.)
Note that the genericity of weak mixing property in the parameter space of translation
flows on translation surfaces is recently proved by A. Avila and G. Forni \cite{Avila}.
A horizontal-vertical polygon is a polygon all of whose sides are in horizontal or
vertical direction, and a horizontal-vertical translation surface is, roughly speaking,
an oriented surface constructed as compact subset on Euclidian 2D space with edges
$\{ l_i^{\pm }\} $ where for each $i$, $l_i^+$ and $l_i^-$ are parallel and have the same
length, and each pair $(l_i^+,l_i^-)$ are topologically connected to each other
and each "jump" from $l_i^-$ to $l_i^+$ is horizontal or vertical in covering space.
Any other results for weak mixing were obtained numerically, see the "square-ring billiard"
\cite{AGR}; see also related numerical studies on rational \cite{A} and irrational
triangles \cite{ACG,CP}. \par
Billiard in a horizontal-vertical polygon whose ratio of length of horizontal or vertical
sides are rational respectively are called to be almost integrable (or A-integrable)
billiard. We can extend this concept for horizontal-vertical translation sufaces;
for horizontal-vertical translation surface, if the ratio of length of horizontal or
vertical "jump" are rational respectively, then we call the translation surface almost
integrable. Dynamics of almost integrable billiards have intensively studied, especially
the dynamics of barrier billiards are well understood, having both discrete and singular
continuous spectrum \cite{R,W} and calculability to fractal dimension (or decay rate of
power-law correlation) of trajectries with quadratic irrational gradients were rigorously
proven \cite{HMcC}, and self-similarity of a trajectry were reported \cite{CO};
see also references therein. \par

In this article, the author report a new example of pseudointegrable systems with weak mixing,
although the system is not horizontal-vertical type. To the best of the author's knowledge,
this is the first time to give an explicit example of pseudointegrable systems with weak
mixing without horizontal-vertical property. \par
The outline of this article is as follows. In Section 2 we describe the concept of translation surface which
appears as invariant surfaces of billiard flows on rational polygons.

The concept of translation surface will be used as more
general meaning in Section 4 and 5, that is, as orientable surface which have finitely many punctures and have flat metric exept
to the punctures. In Section 3 we describe the concept of interval exchange, and then we review Veech's weak mixing theorem for
interval exchange \cite{V1} which gives us a sufficient condition to weak mixing for interval exchanges. In section 4 the main
result is stated, the sufficient condition to have weak mixing property for pseudointegrable systems is stated. 
Here the concept of pseudointegrability is extended to dynamics of geodesic flows on translation surfaces.
To make sense of this sufficient condition, an example whose weak mixing is proved by the condition is proposed in Section 5.
The conclusion is stated in Section 6. In Appendix we review spectral properties of horizontal-vertical billiards.

\section{Translation Surfaces}
Let $P$ be a plane polygon, and let $T_{full}(P)$ be phase space of the billiard inside $P$. Then we can denote
$T_{full}(P) = {P \times {\mathbf R}^2}$, here ${\mathbf R}^2$ represents the space of velocities. Because of motions
of billiards, absolute value of velocity $|v|$ is conservative on time evolution, so $T_{|v|}(P)\subset T_{full}(P)$ are invariant,
where $T_{|v|}(P)$ is defined by fixing the absolute value of velocity to $|v|$, and each
$T_{|v|}(P) (v\neq 0)$ is equivalent to $T_1(P)$, which means $T_{|v|=1}(P)$, apart from time scaling. So we should only treat
the dynamics on $T_1(P)$. From here, $T(P)$ will be used instead of $T_1(P)$ as the phase space. Clearly,
$T(P)=P\times {\mathbf S}^1$, where ${\mathbf S}^1$ stands for set of directions. \par
If $P$ is a rational polygon, any orbit on $T(P)$ has finite directions. For example, there exist four directions in generic orbit
of rectangular or barrier billiards, while there exist six directions in generic orbit of regular triangular billiards.
Each of these numbers of directions is equal to order of Coxseter group $G_c(P)$ of the billiard table $P$.
Here Coxseter group is a discrete group each of whose elements is a rotation or a reflection acting on ${\mathbf S}^1$.
And $G_c(P)$ is defined as the group generated by reflections on the sides of $P$. \par
We can construct the invariant surfaces $S_\theta (P)=P\times\Gamma_P(\theta )\subset T(P)$ for $\theta\in {\mathbf S}^1$ with
$\Gamma_P(\theta )=\{ g\theta :g\in G_c(P)\} $. Naturally, if
$\theta^{\prime }\in\Gamma_P(\theta )$, then $S_{\theta^{\prime }}(P)=S_{\theta }(P)$. In the cases of rectangular or barrier (or
more generally, horizontal-vertical) billiards, for example, $G_c(P)$ is generated by
\begin{equation}
g_h: \theta\mapsto -\theta, \qquad g_v: \theta\mapsto\pi -\theta,
\label{ActionOfReflection}
\end{equation}
and $\# (G_c(P))=\# \{ {\it id.},\ g_v,\ g_h,\ g_v\cdot g_h\} =4$. Here $\# (\cdot )$
represents order of the set. \par
Flows on $S_\theta (P)$ for generic $\theta $ can be thought as geodesic flows on
a Riemannian surface $S(P)$ with flat metric except to finitely many punctures.
Such a surface is called {\em translation surface}.
Let $s_j$ be sides of $P$ and let $g_j\in G_c(P)$ be the reflection on $s_j$. Then
the translation surfaces $S(P)$ related to rational polygons $P$ is defined as follow;
\begin{Definition}{Translation Surfaces of Rational Polygons}
{\it The translation surface $S(P)$ of $P$ is a topological surface
with a metric defined as follows: \\
{\bf 1. Topological sense.} $S(P)$ is topologically equivalent to
$P\times G_c(P)=\{ P(g)\}_{g\in G_c(P)}$, where every $P(g)$ is copy of $P$ and $j$-th
side $s_j(P(g))$ of $P(g)$ and $j$-th side $s_j(P(g_j\cdot g))$ of $P(g_j\cdot g)$ are
topologically identified. \\
{\bf 2. Metric sense.} The metric on $S(P)$ is defined as flat metric except to
the vertices.}
\end{Definition}
For generic $\theta\in {\mathbf S}^1/\Gamma_P$, flow $S_{\theta }(P)$ is equivalent to
geodesic flows on $S(P)$ directed to $\theta $. And geodesic flows on $S(P)$ directed
to generic directions are ergodic. Note that an orbit having less than $\# (G_c(P))$
directions in $T(P)$ may correspond to two or more orbits in $S(P)$, but they are out
of discussions in this article because such orbits are peiodic and not generic. \par

Let $l$ be a segment or a straight loop in a translation surface $S$ and let $\theta $
be a direction transversal to $l$. Then the return map onto $l$ of the translation flow
$S_{\theta }$ is a interval exchange, which will be defined in next section.

\section{Interval Exchange and Veech's theorem}
If a map $T:I\to I$, where $I=[0,a)$, $a>0$, is an one-to-one, onto map and deferential
$dT/dx$ is equal to $1$ for Leb.-a.e. $x\in I$, then $T$ is called an
{\em Interval Exchange}. An interval exchange $T$ with $m$ partitions is denoted
by a positive vector $\lambda\in {\mathbf R_+^m}$ and a permutation $\pi $ on
$\{ 1,\ 2,\ \dots ,m\} $; see \cite{V1} and refelences therein.
Here ${\mathbf R}_+^m=\{\lambda\in {\mathbf R}^m:\lambda_j>0,\ 1\leq j\leq m\} $, and
let $|\lambda |$ be $\sum_{j=1}^m\lambda_j$. An interval exchange $T_{(\lambda ,\pi )}$
is {\em irreducible} if $T[0,\tau )=[0,\tau ),\ \tau >0$, implies $\tau =|\lambda |$.
It is equivalent to say $\pi $ is {\em irreducible}, i.e., that
$\pi\{ 1,\dots ,\ k\} =\{ 1,\dots ,\ k\} $ only for $k=m$. $\Pi_m^0$ denotes the set of
irreducible permutations. For an interval exchange $T_{(\lambda ,\pi )}$, we must recall
some notions for describing Veech's weak mixing theorem.
Let $\sigma_{\pi }$ be the permutation on $\{ 0,\ 1,\ \dots ,m\} $ difined as
\[
\sigma_{\pi }i=\pi^{-1}(\pi i+1)-1\qquad (0\leq i\leq m).
\]
Here $\pi $ is extended to the permutation on $\{ 0,\ 1,\dots ,\ m+1\} $ as
$\pi (0)=0,\ \pi (m+1)=m+1$. Then $\{ 0,\ 1,\ \dots ,m\} $ are decomposed some
$\sigma_{\pi }$-invariant subsets $\Sigma_{\pi }=\{ S_1,\ S_2,\ \dots ,S_r\} $.
There are integral vertors $b_{S_j}$ connected with each $S_j$ by
\[
b_{S_j,i}=\chi_{S_j}(i-1)-\chi_{S_j}(i)\qquad (1\leq i\leq m).
\]
Where $\chi_S$ denotes characteristic function. \par
Adopt the convection that $e(t)=\exp (2\pi it),\ t\in {\mathbf R}^m$. If
$\lambda\in {\mathbf R}_+^m,\ \nu\in {\mathbf R}^m$ set up on
$I^\lambda \equiv [0,|\lambda |)$ a function, $\phi =\phi_{(\lambda ,\nu )}$, as
\[
\phi (x)=e(\nu_j)\qquad (x\in I_j^\lambda ,\quad 1\leq j\leq m)
\]
Here $I_j^\lambda $ is defined as
\[
I_j^\lambda =\left\{
\begin{array}{cr}
[0,\lambda_1) & (j=1), \\
\left[
\sum_{i=1}^{j-1}\lambda_i ,\sum_{i=1}^j\lambda_i
\right) & (\mbox{else}).
\end{array}
\right.
\]
Given also $\pi\in\Pi_m^0$, set $T=T_{(\lambda ,\pi )}$, and suppose $f$ is a measurable,
complex-valued solution to the equation
\begin{equation}
f(Tx)=\phi (x)f(x)\qquad (\mbox{Leb.-a.e.\ }x\in I^\lambda )
\label{Veech}
\end{equation}
Then,
\begin{Theorem}{Veech's weak mixing theorem \cite{V1}}
{\em Let $\pi\in\Pi_m^0$. For a.e. $\lambda\in {\mathbf R}_+^m$ it is true for all
$\nu\in {\mathbf R}^m$ that if} (\ref{Veech}) {\em admits a nonzero measurable
solution, then $b_S\cdot\nu\in {\mathbf Z}$ for all $S\in\Sigma (\pi )$.}
\label{Veech'sWeakMixing}
\end{Theorem}
Or conversely,
\begin{Proposition}
{\em Let $\pi\in\Pi_m^0$. For a.e. $\lambda\in {\mathbf R}_+^m$ it is true for all
$\nu\in {\mathbf R}^m$ that if $S\in\Sigma (\pi )$ exists such that
$b_S\cdot\nu\notin {\mathbf Z}$, then} (\ref{Veech}) {\em does not admit any nonzero
measurable solution.}
\end{Proposition}

With $\nu =(1,\ 1,\dots ,1)$, this sufficient condition for weak mixing works on interval exchanges. \par
In next section, how this theorem works on pseudointegrable systems is described.

\section{A New Approach to Weak Mixing for PseudoIntegrable Systems}
Let $S_\theta $ be an ergodic component of a translation surface $S$, and $l$ be a
segment or a straight loop on $S$ transversal to $\theta $. Then the return map of
$S_\theta $ onto $l$ is an interval exchange $T_{(\lambda ,\pi )}$ with $m$
disconnected points, where $m$ is the number of punctures of $S$. It is obvious from
the construction of $T_{(\lambda ,\pi )}$ that the length of $l$ is equal to
$|\lambda |$, and for each connecting component $l_i\subset l$, any points
$x,y\in l_i$ are contemporaneously returned onto $l$ by flows on $S_\theta $. So the
characteristic value problem of unitary operator $U^t:L^2(S)\to L^2(S)$
\begin{equation}
U^tf(x)=\exp (2\pi i\alpha )f(x),\qquad f\in L^2(S),
\label{continuousEigen}
\end{equation}
on the translation flow $S_\theta $ on $S$ is reduced on the interval exchange
$T_{(\theta ,\pi )}$ on $l$ as the form of
\begin{equation}
\widehat{U}f_l(x)=\phi (x)f_l(x),\qquad f_l\in L^2(l).
\label{discreteEigen}
\end{equation}
In the equation (\ref{continuousEigen}) and (\ref{discreteEigen}), $U^t$, $\widehat{U}$
are defined as $U^tf(x)=f(\varphi^tx)$, $f\in L^2(S)$ where $\varphi $ is the time
evolution by the flow of $S_\theta $ and as $\widehat{U}f(x)=f(T_{(\lambda ,\pi )}x)$,
$f\in L^2(l)$ respectively, $\phi (x)$ in (\ref{discreteEigen}) corresponds to $\alpha $
in (\ref{continuousEigen}) as
\begin{equation}
\phi (x)=e(\alpha t_j),\qquad (x\in I_j^\lambda\quad 1\leq j\leq m),
\label{phi}
\end{equation}
and $f_l$ in (\ref{discreteEigen}) is defined from $f$ in (\ref{continuousEigen}) by
restricting onto $l$;
\[
f_l(x)=f(x),\qquad (x\in l\subset S).
\]
To compare (\ref{Veech}) and (\ref{discreteEigen}), we can observe that the equation
(\ref{discreteEigen}) is of the form (\ref{Veech}).
Therefore Veech's weak mixing theorem works on $S_\theta $ in the sense describing bellow. \par
For ergodic $S_\theta $, the map $T_{(\lambda ,\pi )}$ on any $l$ is naturally ergodic and,
therefore, the permutation $\pi $ is irreducible. So, it is true from Veech's weak mixing
theorem that
\begin{Corroraly}
{\em For $T_{(\lambda ,\pi )}$ on $l$ induced from an ergodic $S_\theta $, if
$S\in\Sigma_\pi $ satisfies the condition $b_S\cdot\alpha t\notin {\mathbf Z}$,
where $t$ is the positive vector whose $j$-th component is defined as $t_j$ in}
(\ref{phi}), {\em then} (\ref{continuousEigen}) {\em for $S_\theta $ doesn't
have the characteristic value $\alpha $ and, moreover, doesn't have the $\alpha n$
for any $n\in {\mathbf Z}$.}
\end{Corroraly}
This collorary is available only if $b_S\neq 0$ exists for some $S\in\Sigma_\pi $. \par
If $S_j,S_k\in\Sigma_\pi $ exist such that $b_{S_j}\cdot t$ and $b_{S_k}\cdot t$ are
linearly independent with integral coefficient, then the value $\beta $ such that
$b_{S_j}\cdot t\beta\in {\mathbf Z}$ isn't characteristic value because of that
$b_{S_k}\cdot t\beta\notin {\mathbf Z}$. So, such an ergodic component is weakly
mixing. \par
For the last of this section, let us summarize the above condition to weak mixing;
\begin{Theorem}{A Sufficient Condition for Weak Mixing}
{\em If $S_j,S_k\in\Sigma_\pi $ exist such that $b_{S_j}\cdot t$ and
$b_{S_k}\cdot t$ are linearly independent with integral coefficients, then the flow
$S_\theta $ is weak mixing on the ergodic component.}
\label{MainTheorem}
\end{Theorem}
To indicate the availavility of this condition we will show an example of weak mixing
pseudointegrable system whose weak mixing property is supported from above condition
of weak mixing.

\section{Example}
In this section, we show an algorithm to construct pseudointegrable systems with weak
mixing based on the sufficient condition in the last section. Note that this algorithm
is incomplete to construct such systems in many points. \par
To construct such a system,
first, we must find such a permutation $\pi $ that there exist $b_{S_j}$ and $b_{S_k}$
($S_j, S_k\in\Sigma_\pi ,\ j\neq k$) and they are linearly independent.
Since $\sum_ib_{S_i}=0$, $\# (\Sigma_\pi )$ must be $\geq 3$. As such an example
of $\pi $, let us choose
\begin{equation}
\pi =\left(
\begin{array}{c}
1\ \ 2\ \ 3\ \ 4 \\
4\ \ 2\ \ 3\ \ 1
\end{array}
\right) .
\label{permutation}
\end{equation}
Next, we must find a corresponding translation surface $S$, a direction of flow $\theta $
and a segment (or the loop) $l$, of which permutation $\pi $ of the
interval exchange $T_{(\lambda ,\pi )}$ will become the permutation finded above.
For our example of $S_\theta $ and $l$ corresponding to our permutation (\ref{permutation}),
let us choose them as describing in
\begin{figure}
\begin{center}
\includegraphics[height=5cm]{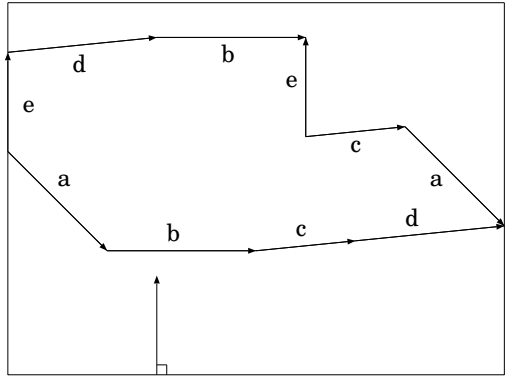}
\caption{A translation surface with weak mixing for $\theta = \pi /2$ and for the other
parameters as generic. Each pair of vectors $a,\ b,\dots ,e$ are identified each other.
This system is parametrized by the 2D vectors $a,\ b,\dots ,e$. And the outer square is
unit square and horizontal or vertical sides of this square are identified respectively.}
\label{FigA}
\end{center}
\end{figure}
Fig.\ref{FigA} and as $x$-axis. \par
The last, we must check whether there exist $b_{S_j}\cdot t$ and $b_{S_k}\cdot t$
($S_j, S_k\in\Sigma_\pi ,\ j\neq k$) such that they are linearly independent with integral
coefficients. If they exist, then it is supposed by theorem \ref{MainTheorem} that
the flow $S_\theta $ is weak mixing. For our case, weak mixing is immediately proven for
generic parameters. More explicitly, $S_\theta $ is ergodic if horizontal components of
$\{ a,\ b,\dots ,e\} $ are linearly independent with integral coefficients, and in addidion
if vertical component of $e$ is irrational, then $S_\theta $ is weak mixing.

\section{Conclusion}
On pseudointegrable systems, we were able to show weak mixing property
only for horizontal-vertical cases. On the other hand,
Veech's weak mixing theorem for interval exhange was used only for 
A-integrable systems to prove that the complement spectra of
immediately given discrete component is singular continuous \cite{ZE}. \par
In this paper, the author established that Veech's weak mixing theorem
(Theorem \ref{Veech'sWeakMixing}) can be applicable to pseudointegrable
systems for proving their weak mixing property, and that is stated
as Theorem \ref{MainTheorem}
To make sense of theorem \ref{MainTheorem}, it is proposed that there
really exists a translation flow of a translation surface
whose weak mixing property is supported by this theorem.
Furthermore, this translation surface is not horizontal-vertical type.
This approach to weak mixing for pseudointegrable systems is the defferent
way from previous works in this field. \par
The example which is proposed in section 5 is not billiard system,
so it is still opened whether there exist some examples whose weak mixing
property is supported by this theorem. \par
Note that this aproach has further problems. First, theorem \ref{MainTheorem}
gives us only the sufficient condition for weak mixing property.
Second, the algorithm shown in Section 5 is incomplete too, in other words,
even if the algorithm failed with some segment $l$ in $S$,
the algorithm may succeed with other segment $l$. \par
To improve the first incompleteness, we have to extend Veech's weak mixing
theorem (theorem \ref{Veech'sWeakMixing}) to sufficient and necessary
condition. And to improve the second incompleteness, we have to study
the relation between interval exchange transformations on different segments
in the same translation surface.

\section{Acknowledgement}
The Author would like to thank Tetsuro Konishi for checking this paper and giving many helpful
advices for the author, and to thank Masashi Tachikawa for recommending this subject
to the author.

\appendix
\section*{Appendix: Spectral Properties of Horizontal-Vertical Translation Surfaces}
\begin{figure}
\begin{center}
\includegraphics[height=5cm]{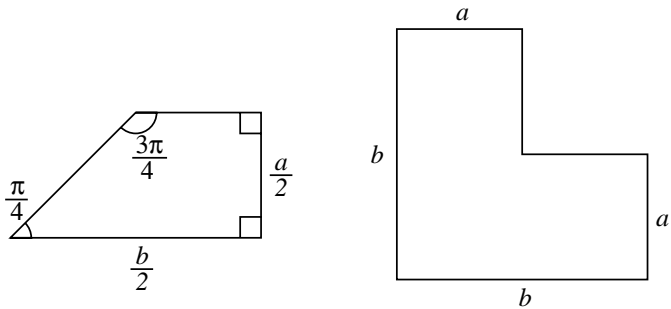}
\caption{A horizontal-vertical billiard table (left)
and its translation surface (right) are presented. In the translation surface,
every point on the horizontal, or vertical, sides are topologically identified to
the opposite point which has the same value of $x$-, or $y$-, component,
respectively. For generic $\theta $, if $a/b$ is rational $S_\theta $ is
almost integrable, while if $a/b$ is irrational $S_\theta $ is weak mixing.}
\label{FigB}
\end{center}
\end{figure}

In this section, we show the spectral properties of horizontal-vertical
translation surfaces, whose definition is already stated in Section 2.
We discuss only about the translation surfaces described in Fig.\ref{FigB},
but any other cases can be discussed similarly. \par
Let $S$ be the translation surface described in Fig.\ref{FigB} and $\theta $ be
the direction of which $S_\theta $ is ergodic flow. Then the flow $S_\theta $ is
almost integrable if and only if $a/b$ is rational. \par
As the almost integrable case, let us choose the parameters $a,\ b$ in Fig.\ref{FigB}
as $b=2a$. The vector field on $S_\theta $ is $v=(\cos\theta ,\ \sin\theta )$ uniformly.
Therefore the eigenequation (\ref{continuousEigen}) admits infinitly many eigenfunctions
$\{f_{jk}\}_{j,k\in {\mathbf Z}}$
\begin{equation}
f_{0;jk}(x,y)=e(jx+ky),
\label{eigenfunction}
\end{equation}
with corresponding eigenvalues
\begin{equation}
\alpha_{jk}=j\cos\theta +k\sin\theta .
\label{eigenvalues}
\end{equation}
Here $e(t)=e^{2\pi it}$, which appeared in section 3.
Any other eigenfunction couldn't be admitted by the equation (\ref{continuousEigen}).
To prove this, let $L_0^2\subset L^2(S)$ be the subspace spanned by the eigenfunctions
(\ref{eigenfunction}). Then the orthogonal complement subspace $L_0^{2\perp }$ is
decomposed into two subspaces $L_1^2$ and $L_2^2$ which are orthogonal to each other.
Where $L_l^2,\ l=1,2$, are spanned by
\[
f_{l;jk}(x,y)=\left\{
\begin{array}{c}
e(jx+ky)\qquad (0\leq x,y<1) \\
e(jx+ky+l/3)\qquad (0\leq y<1\leq x<2) \\
e(jx+ky+2l/3)\qquad (0\leq x<1\leq y<2)
\end{array}
\right.
\]
respectively. Suppose $f\in\mbox{span}\{ f_{l;jk}:l=1,2,\ j,k\in {\mathbf Z}\} $ and
$U^tf=e(\alpha t)f$, then
\begin{equation}
|f|=\mbox{const.}\qquad \mbox{Leb.-a.e.}\ (x,y)\in S
\label{const}
\end{equation}
is obvious from ergodicity, and
\begin{equation}
f(x,y)+f(x+1,y)+f(x,y+1)=0,\qquad\mbox{Leb.-a.e.}\ (x,y)\in [0,1)\times [0,1)
\label{zero}
\end{equation}
is immediately obtained by the definition of $f_{l;jk}$. Moreover $f$ have to be
$\in L_1^2$ or $\in L_2^2$ because if there exist $f_1\in L_1^2$ and $f_2\in L_2^2$
($||f_1||, ||f_2||\neq 0$) such that $f=f_1+f_2$, then
\[
f^2=2f_1f_2+f_2^2+f_1^2
\]
is also eigenfunction, however, $2f_1f_2\in L_0^2$ and $f_2^2+f_1^2\in
L_0^{2\perp }$, so both $2f_1f_2$ and $f_2^2+f_1^2$ are the eigenfunctions having
the same eigenvalue because of the invariance of $L_0^2$ under the action $U$,
this situation is impossible because of erogodicity.
So, we suppose that $f\in L_1^2\ (\mbox{mod}\ 0)$ (or $\in L_2^2\ (\mbox{mod}\ 0)$).
Its orbit $U^tf$ have to be $\in L_1^2$ ($\in L_2^2$, respectively) in any
$t\in {\mathbf R}$.
For any point $p_0=(x_0,y_0)\in A_0=\{ (x,y)\in (0,1)\times (0,1)\ |
\ y-1-(x-1)\tan\theta <0\} $,
we can take the positive numbers $t_0>0,\ \varepsilon >0$ such that
${\mathcal U}_\varepsilon (p_0)\subset (0,1)\times (0,1)$, where
${\mathcal U}_\varepsilon $ is the $\varepsilon $-neighberhood of $p_0$, and
$S_\theta ^t{\mathcal U}_\varepsilon (p_0)\subset (0,2)\times (0,1)$ for any
$0\leq t\leq t_0$, and $S_\theta ^{t_0}{\mathcal U}_\varepsilon (p_0)\subset
(1,2)\times (0,1)$. Then, for $p_1=(x_0+1,y_0)$ and $p_2=(x_0,y_0+1)$ it is
clear that ${\mathcal U}_\varepsilon (p_1)\subset (1,2)\times (0,1),
\ {\mathcal U}_\varepsilon (p_2)\subset (0,1)\times (1,2),
\ S_\theta ^{t_0}{\mathcal U}_\varepsilon (p_1)\subset (0,1)\times (0,1)$ and
$S_\theta ^{t_0}{\mathcal U}_\varepsilon (p_2)\subset (0,1)\times (1,2)$.
So we can conclude that $f|_{\mathcal U}$, here
${\mathcal U}={\mathcal U}_\varepsilon (p_0)\cup{\mathcal U}_\varepsilon (p_1)
\cup{\mathcal U}_\varepsilon (p_2)$, have to be equal to $0$, thus
$f|_{A_0}=0\ (\mbox{mod}\ 0)$.
We may conclude the same result for
$f|_{(0,1)\times (0,1)-A_0}$ by the same way. So, $f|_{(0,1)\times (0,1)}=0
\ (\mbox{mod}\ 0)$, and thus, $f=0\ (\mbox{mod}\ 0)$,
but such a function couldn't be an eignfunction.
Therefore, we obtain the result that any other eigenfunction does no exist. \par
In the case of that $a/b$ is irrational, there is the sequence
$\{ n_i/m_i\}_{i\in {\mathbf N}}$ which is the rational approximation by continued
fraction. In each case of $a/b=n_i/m_i$, the system is A-integrable and there is
the subspace $L_{d.i}^2(S)\subset L^2(S)$ which is spanned by the eigenfunctions.
In each $n_i/m_i$, $S$ is devided to $N_i$ copies of fundamental rectangle.
(In our case, $N_i=m_i^2-n_i^2$, and the fundamental rectangle is sized as
$\frac{b}{m_i}\times\frac{b}{m_i}$.) Then,
\[
N_i\to\infty ,\qquad (i\to\infty ),
\]
and $dim(L_{d,i}^2(S))$ shrinks more and more with $i$ getting large.
As $i\to\infty $, $dim(L_{d,i}^2(S))\to 0$. \par
So we reach the conclusion that an ergodic flow $S_\theta $ is A-integrable
or weak mixing for any horizontal-vertical translation surface $S$.

\end{document}